\documentclass{article}
\usepackage[utf8]{inputenc}
\usepackage{arxiv}
\usepackage[utf8]{inputenc} 
\usepackage[T1]{fontenc}    
\usepackage{amssymb}
\usepackage{amsthm}
\usepackage{lineno}
\usepackage{algorithm}
\usepackage{algorithmic}
\usepackage{graphicx}
\graphicspath{{./images/}}
\usepackage{amsmath}
\usepackage{amsfonts}
\usepackage{mathptmx}
\usepackage{multirow}
\usepackage{mathrsfs}
\usepackage{colortbl}
\usepackage{subfigure}
\usepackage{ulem}

\newtheorem{remark}{Remark}[section]

\title{A Deep Neural Network for Knot Placement in B-spline Approximation}

\author{Jiaqi Luo \\
	School of Mathematical Sciences\\
	Soochow University\\
	No.1 Shizi Street, Suzhou, Jiangsu Province,  China\\
	Data Science Research Center\\
	Duke Kunshan University\\
	No.8 Duke Ave, Kunshan, Jiangsu Province, China\\
	\texttt{jiaqi.luo@dukekunshan.edu.cn} \\
	\And 
	Zepeng Wen\\
	School of Mathematical Sciences\\
	Soochow University\\
	No.1 Shizi Street, Suzhou, Jiangsu Province,  China\\
	\texttt{zpwen@stu.suda.edu.cn}
	\And
	Hongmei Kang\\
	School of Mathematical Sciences\\
	Soochow University\\
	No.1 Shizi Street, Suzhou, Jiangsu Province,  China\\
	\texttt{khm@suda.edu.cn} \\
	\And 
	Zhouwang Yang\\
	School of Mathematical Sciences\\
	University of Science and Technology of China\\
	No.96 Jinzhai Road, Hefei, Anhui Province, China\\
	\texttt{yangzw@ustc.edu.cn} \\
	 }
\date{May 2022}

\begin{document}
\maketitle
\begin{abstract}
Automatically determining knot number and positions is a fundamental and challenging problem in B-spline approximation.
In this paper, the knot placement is abstracted as a mapping from initial knots to the optimal knots. We innovatively introduce a deep neural network solver to approximate the mapping. The neural network is composed of several subnetworks. 
Each subnetwork is designed to approximate the optimal knot positions in the case of fixed knot number. All the subnetworks are stacked together to find the optimal knots (including knot number and knot positions) within some given tolerance. Owing to the powerful approximation capabilities, as well as mature algorithms developed in deep learning, the proposed method can effectively and efficiently find the optimal knot number and knot positions. Several numerical examples are demonstrated to show the superiority of our method.
\end{abstract}



\section{Introduction}
A B-spline curve $c(t)$ of degree $p$ defined over a knot vector $\textbf{\textit{U}}=\{t_0,t_1$, $\cdots,t_{n+p+1}\}$ is defined as 
\begin{equation}
\label{E:bsplinecurve}
\textbf{C}(t)=\sum_{i=0}^n\textbf{\textit{C}}_i~N_i^p(t),
\end{equation}
where $\textbf{\textit{C}}=\{\textbf{\textit{C}}_0,\textbf{\textit{C}}_1,\cdots, \textbf{\textit{C}}_n\}$ be $n+1$ control points in $\mathbb{R}^d$ and $N_i^p(t)$ be the $i$-th B-spline basis functions
of degree $p$. In B-spline fitting, the degree $p$ is most often fixed as $p=2$ or $p=3$. Therefore, knots and control points are two components to determine the B-spline that satisfies fitting requirements. Once $U$ is fixed, B-spline fitting degenerates to a linear least-squares problem, 
\begin{equation}
\label{E:lseq}
\min_{\textbf{\textit{C}}}\|\textbf{\textit{P}}-\textbf{\textit{AC}}\|_F,
\end{equation}
where $\textbf{\textit{P}}=(P_0, P_1, ... , P_N)^T$ is a given data point vector,   $\textbf{\textbf{A}}=(a_{ij})_{N \times (n+1)}$ with $a_{ij}=N_j^p(s_i)$ is the coefficient matrix, $s_i$ is the corresponding parametric value of $P_i$,   
$\textbf{\textit{C}}=(C_0,C_1,...,C_n)^T$ is the control point vector, and $||\cdot||_F$ denotes the Frobenius norm. 

Obviously, the adequate approximations by \eqref{E:lseq} are built on appropriately chosen knot vectors. Uniform knots or artificially specified knots are far from sufficient to achieve the goal. This can also be confirmed from the numerical experiments shown in Section 4. The fact that free knots are vital to improve fitting performance has been recongnized in the early works \cite{jupp1978approximation,burchard1974splines,de1968least2}. The B-spline fitting problem with free(variable) knots is generally formulated as,
\begin{equation}
\label{E:topprob}
\min_{n,\textbf{\textit{C}},\textbf{\textit{U}}}\|\textbf{\textit{P}}-\textbf{\textit{A(U)C}}\|.
\end{equation}
The knot number $n$, knot positions $\textbf{\textit{U}}$ and control points $\textbf{\textit{C}}$ are all unknowns. Alternatively, the more practical version is,
\begin{align}
	\label{E:inv1}
	&~~\min_{\textbf{\textit{U, C}}}~n\\ 
	&s.t.~~ \|\textbf{\textit{P}} - \textbf{\textit{A(U)C}}\|  \leq \epsilon, \notag
\end{align}
where $\epsilon$ is the fitting tolerance required in practice. Because of the nonlinear relationship between the knots of B-splines, both the formulations  \eqref{E:topprob} and \eqref{E:inv1} are non-convex and highly nonlinear optimization problem, and the numerical solutions are very hard to obtain.


In the existing literature, the following two easy problems are mostly considered: 
\begin{itemize}
	\item[\textbf{P1.}] finding the approximation which has the minimum fitting error in the case of fixed knot number, which is formulated as
	\begin{equation}
	\label{E:lseq2}
	\min_{\textbf{\textit{C}},\textbf{\textit{U}}}\|\textbf{\textit{P}}-\textbf{\textit{A(U)C}}\|.
	\end{equation}
	\item[\textbf{P2.}] finding the approximation with  as few knots as possible under a certain error tolerance $\epsilon$.

\end{itemize}
The first problem (\textbf{P1}) is often solved by means of nonlinear optimization techniques \cite{jupp1978approximation,globalopt,yoshimoto2003data,ilker2009automatic}. The good performance of the approximations is based on a good guess about the knot number and initial knot positions, which are not easy to give in practice. 
While the second one (\textbf{P2})  is generally solved by heuristic methods, such as adaptively inserting knots \cite{park2007b-spline}, adaptively removing knots \cite{lyche1987knot,lyche1988a} or selecting from given initial knots \cite{yuan2013adaptive}. The knots calculated by heuristic methods are usually redundant than expected.


In this paper, we propose a deep learning method to address the knot placement problem. The proposed deep neural network is composed of several subnetworks. Each subnetwork is designed to approximate the optimal knot positions in the case of fixed knot number. All the subnetworks are stacked together to find the optimal knot number and knot positions within some given tolerance. The  loss function of the entire neural network is defined by the fitting error. With the proposed neural network, the knot number and knot positions are determined simultaneously, as well as the corresponding B-spline approximations meet the error tolerance.


We give an intuitive explanation of the key idea of the proposed method. On the one hand, the problem (\textbf{P1}) is innovatively understood as a mapping from initial knots to optimal knots. As you can image, this mapping is complicated, and it is a good idea to approximate it through deep neural networks, since deep neural networks can be understood as a composite function of nonlinear transformations. All the parameters of nonlinear transformations constitute the parameters of the neural network. Generally the parameters of neural networks are far more than the DOFs (degrees of freedom) in problem \eqref{E:inv1}. This over-parametrized nature makes the deep neural networks more likely to find the global optimum, which has been explored and analysized in many papers \cite{deeptheory1,deeptheory2,deeptheory3}.
Moreover, the well-developed optimization algorithm in deep learning accelerates the solving process. The proposed method is much faster than the current existing methods.
On the other hand, a deep neural network can only give the optimal knot positions for fixed knot number,  so we stack several sub-neural networs together to assess the effect of different knot number on the fitting error.

\section{Related Work}
The nonlinear optimization problem   (\textbf{P1})  usually has multiple local minima, thus global optimum is very hard to achieve.  In the earlier work \cite{jupp1978approximation}, the alternative iteration method is proposed for solving the optimization problem. The modern techniques developed in nonlinear optimization field, such as the cutting angle method of
deterministic global optimisation \cite{globalopt}, the particle swarm optimization method \cite{galvez2011efficient,galvez2012particle}, the genetic
algorithm \cite{yoshimoto2003data}, and the artificial immune system algorithm \cite{ilker2009automatic}, are adopted in solving (\textbf{P1}) .

Many reseachers have made a lot of effort to solve the problem (\textbf{P2}). In \cite{lyche1987knot,lyche1988a}, the authors remove knots adaptively based on the assigned weights, which are used to measure the significance of the knots in reducing fitting error. Inserting knots adaptively is also a common way to solve the problem (\textbf{P2}). The intuitive idea is that few points should be placed where the curve looks flat, but many where it has sharp features. Discrete curvature is often utilized to characterize this change of curves \cite{crampin1985linear,park2007b-spline}. The discrete curvature of the points is smoothed to expose the characteristics for dense and noise data in \cite{li2005adaptive}. In \cite{ravari2016reconstruction}, the authors use group testing to find salient points and then add a new salient point to the fitting process iteratively until the Akaike Information Criterion (AIC) is met. Laube et al. \cite{laube2018learnt} employ the machine learning method, support vector machine, to choose knots in curve fitting. The curvature information is employed as geometric features therein. In \cite{zhao2011adaptive}, the initial values of the knots are generated by the Monte Carlo method and the optimal knots are searched among the initial candidates based on the fitness. The sparse optimization models are employed to select knots from given knots in \cite{yuan2013adaptive,optimalmin15}. The two-stage frameworks are proposed in \cite{kang2015knot,luo2019knot, luo2021knot}, where initial knot candidates are determined by solving a linear sparse model and the final knots are obtained by adjusting the intial candidates in the second step.


The paper is organized as follows.
The proposed neural network solver for B-spline approximation is presented in detail in Section \ref{s.method}. Several numerical experiments are demonstrated in Section \ref{s.numexp} to show the supriority of the proposed method. Finnaly, we conclude the paper with a discussion of future research problems in Section \ref{s.conclusion}.

%
%

\section{Methodlogy \label{s.method}}
The knot placement task is composed of determining the optimal knot number $n$ and calculating the optimal positions of $n$ knots. We can image that there exists a mapping $g$ from any given initial knot vector $U_0$ to the optimal knot vector $U^*$ (the solution of the problem \eqref{E:topprob} or  \eqref{E:inv1}, that is $g(U_0)=U^*$. Such a mapping is certainly complicated and traditional techniques can not figure out it. The deep neural network (DNN) is a  preferred choice for approximating the mapping $g$. The two vectors $U_0$  and $U^*$ are regarded as the input and output of the network. However, the dimension of the out of DNN is usually fixed, this means one DNN can only give the optimal knots in the case of fixed knot number. Generally the knot number is hard to guess. But the upper and lower bound is easy to give. So we use several deep neural networks to represent the mappings with regard to various knot number, and the fitting error is used for loss function to assess the effect of knot number on fitting performance. The optimal number is determined automatically by the loss function of the deep neural network.   


\subsection{Feedforward neural network model (FNN)}
\begin{figure}[!htbp]
	\centering
	\includegraphics[scale=0.9]{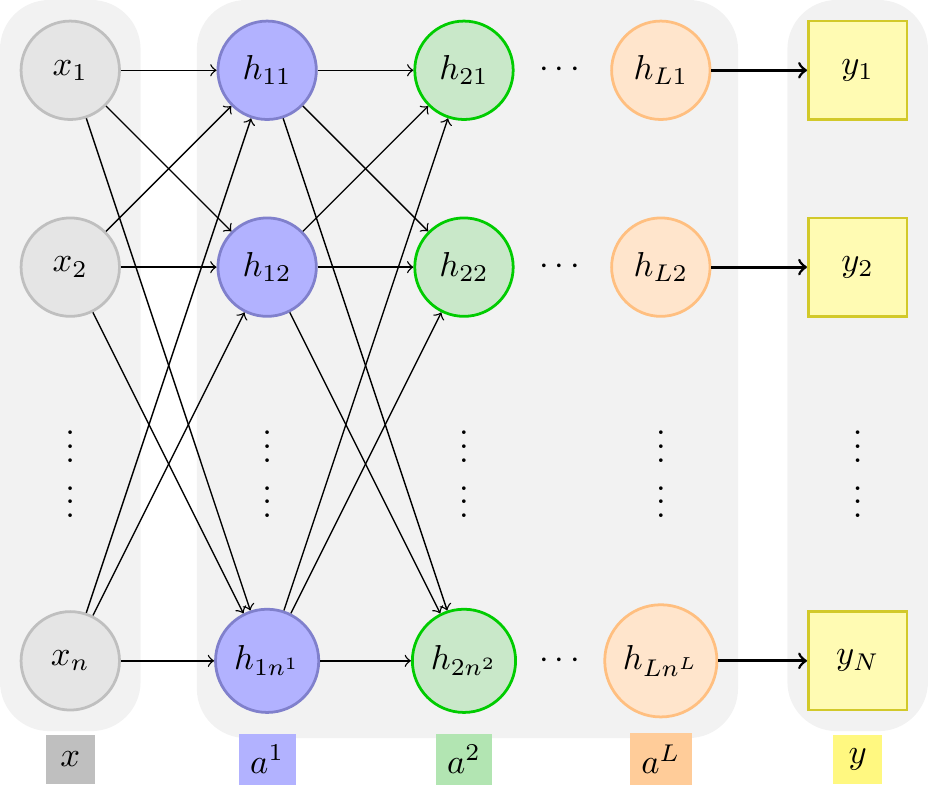}
	\caption{\label{F:dnn} A general feedforward neural network.}
\end{figure}

A FNN model is composed of multiple layers and aims to learn features from input data.
Fig.~\ref{F:dnn} shows the architecture of a general feedforward neural network model with $L$ layers, where $\mathbf{x}=\{x_1,x_2,\cdots,x_n\}$ is the input vector, $\mathbf{a}^l=\{h_{l1},h_{l2},\cdots,h_{ln^l}\}$ is the variable vector on the $j$-th hidden layer, and $\mathbf{y}=\{y_1,y_2,\cdots,y_N\}$ is the output.
The output $\mathbf{y}$ can be understood as a composite function $F(\mathbf{x})$ of several nonlinear functions,
\begin{equation}
	\label{E:mlp}
	\mathbf{y} := F(\mathbf{x}) = (f_{1}^{W_1,b_1} \circ f_{2}^{W_2,b_2} \circ ... \circ f_{L}^{W_L,b_L})(X).
\end{equation}
where $X \in \mathbb{R}^n$ is \textbf{input}, $\mathbf{y} \in \mathbb{R}^N$ is \textbf{output},   $f_{l}^{W_l,b_l} := f_{l}(\mathbf{a}^l)$, $f_l$ is the activation function from the $l$-th layer, $\mathbf{a}^l=W_{l}\mathbf{x}_{l} + b_{l}$, $W_{l}$ and $b_l$ are unknown matrices and vectors. The unknowns $W_{l}$ and $b_l$, $l=1,2,\cdots,L$ are called parameters of a neural network and can be learned by training.

\subsection{Network structure}
\begin{figure*}[!htbp]
	\centering
	\subfigure{\includegraphics[scale=0.4]{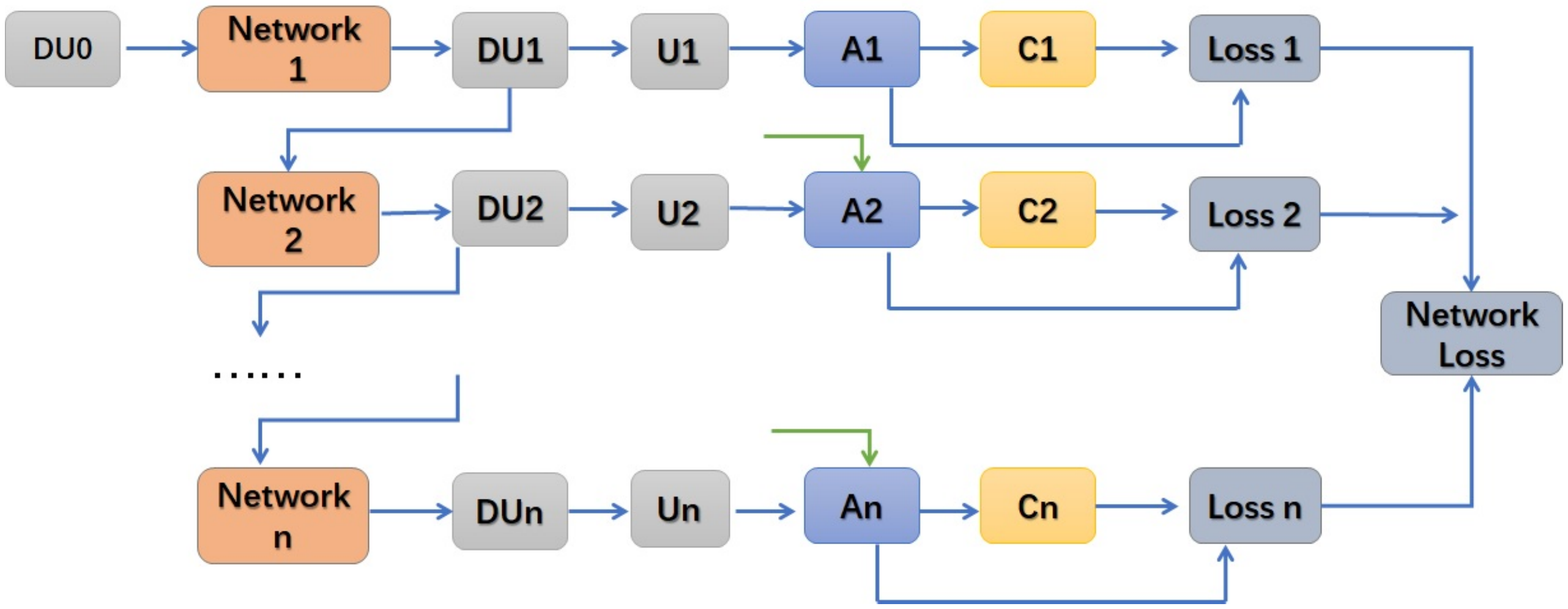}}
	\caption{\label{F:structure} The deep neural network architecture proposed for knot placement.}
\end{figure*}

The construction of an appropriate network structure usually requires a lot of experiments. A too deep neural network would result in long training time and the overfitting problem.  After a lot of experiments, we found that the most simple FNN is enough for the present problem.
Fig.~\ref{F:structure}  illustrates the network architecture used in this paper. The neural network is composed of $n$ subnetworks and each subnetwork has many connections:
\begin{itemize}
\item The network $network_k$, $k=1,2,\cdots, n$ is a multilayer feedforward neural network with two hidden layers. The activation function is chosen as ReLU,
$$ReLU(x) = \max(0, x),$$
except the last one $f_n$, which is set as the Softmax function ($Softmax(x_i) = \frac{\exp(x_i)}{\sum_j{\exp(x_j)}}$). This is done to maintain the monotonicity of the knot sequence. 

\item The input of $network_1$ is the knot difference vector, i.e. $DU_0=\{\bigtriangleup_1, \bigtriangleup_2,\cdots, \bigtriangleup_{n_1}\}^T$, where $\bigtriangleup_j = t_{j+1}-t_j$, $t_j$ are the initial knots given for $network_1$.
The $DU_k$ is the output of the previous subnetwork $network_{k}$ and the input of the next $network_{k+1}$. 

\item In the connection $DU_k \to network_{k+1} \to DU_{k+1}$, $k=1,2,\cdots,n-1$, the dimensions of $DU_k$ and $DU_{k+1}$ are set monotonically decreasing, $dim(DU_{k+1})=dim(DU_{0})-C$ where $C$ is a constant. For example, if $dim(DU_0) = dim(DU_1) = 10$ and $C=1$, then $dim(DU_{k+1}) = 10-k, k \geq 1$. By this  setting, we get a sequence of knot vectors  with different dimensions, $U_i$, $k=1,2,\cdots,n$ .

\item $DU_k \to U_k$ is to recover knots $U_k$ from the difference knot vector $DU_k$. This can be implemented by 
$$U_k = M_k\cdot DU_k,$$
where $M_k$ is a lower triangular matrix with all elements 1.  Notice this connection has no parameters to be optimized.

\item $U_k \to A_k$ is to construct the matrix $A_k=(a_{ij})$ (the same as the one in \eqref{E:lseq}), $a_{ij}=N_i^p(s_j)$, where $N^p_i$ is the B-spline basis function associated with $U_k$. There are also no parameters needed to be optimized.
	
\item $A_k \to C_k$ is to compute the control point vector, $C_k = (A_k^TA_k)^-1A_k^TP$. This connection is only a matrix calculation process and no parameters need to be optimized.
	
\item $C_k, A_k \to Loss_k$ is to define the loss function of $network_k$. The evaluations at $s_j$ of the predicted curve by  $network_k$ are defined as $A_k*C_k$. The loss function is thus defined by $Loss_k = \| P-A_kC_k\|^2_F$.

\item The $Network~Loss$ is defined as the sum of all subnetworks, $Network~Loss = \sum_{k=1}^{n} Loss_k$. 
\end{itemize}

The proposed network gives a sequence of optimal knot vectors, $U_1, U_2,\cdots, U_n$, the knot number of which is gradually reduced. Intuitively, $dim(U_1)$ and $dim(U_n)$ are exactly the guess of lower and upper bounds of the optimal knot number. The optimal knot number is determined by the loss function curves. In order to give more compact bounds, a rough train of the proposed network is suggested. In a rough train phase, the initial knots $U_0$ (corresponding to the knot difference vector $DU_0$) are given uniformly and $dim(U_0)$ is set to be a larger number.
The  two hyper-parameters $C$ and $n$ are set to be $C=5$ and $n=6$ in our experiments.

\begin{figure*}[!htbp]
\centering
\subfigure[Given data and B-spline approximation]{\includegraphics[scale=0.3]{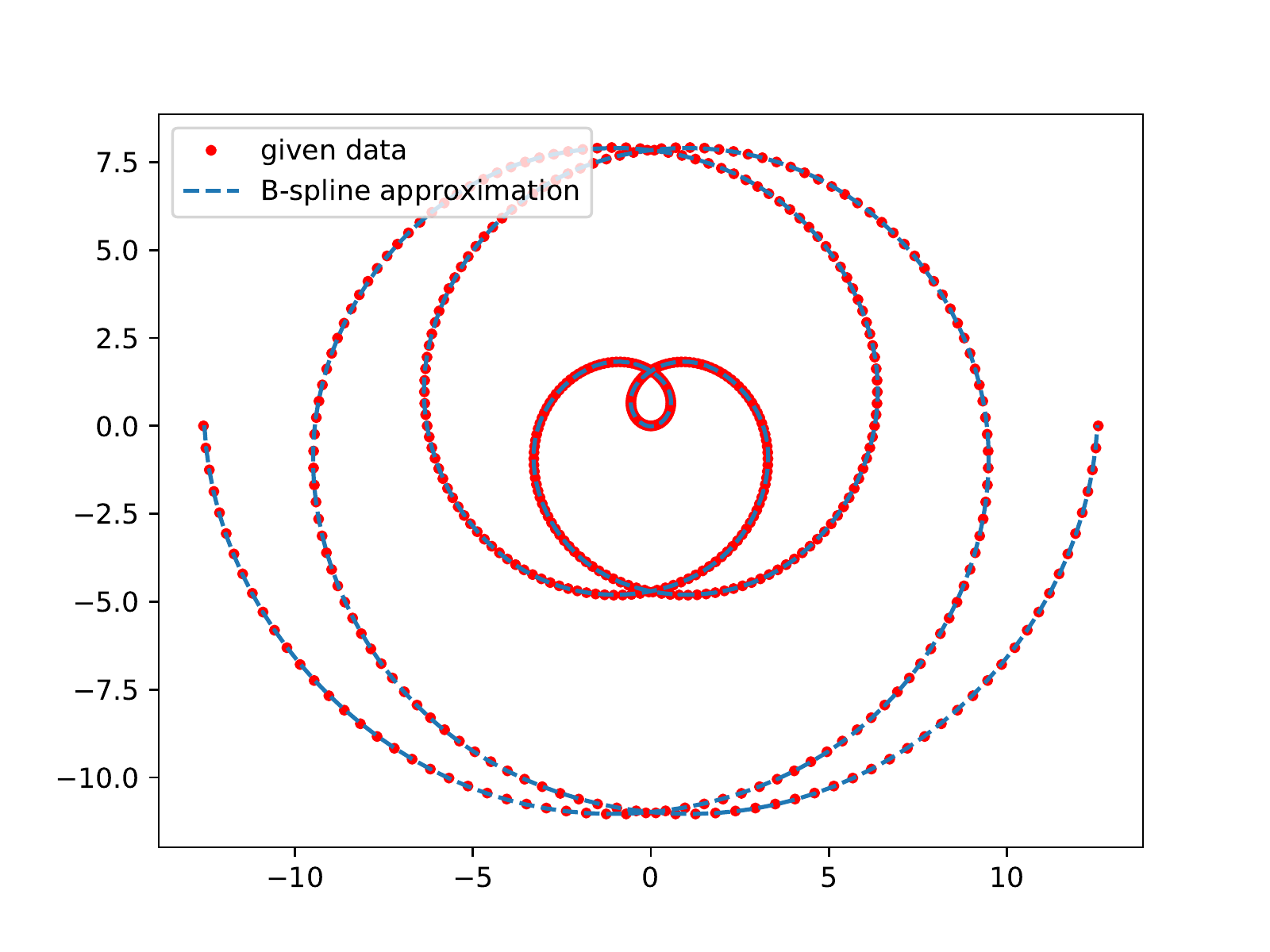}}
\subfigure[Loss curves in the rough train phase]{\includegraphics[scale=0.3]{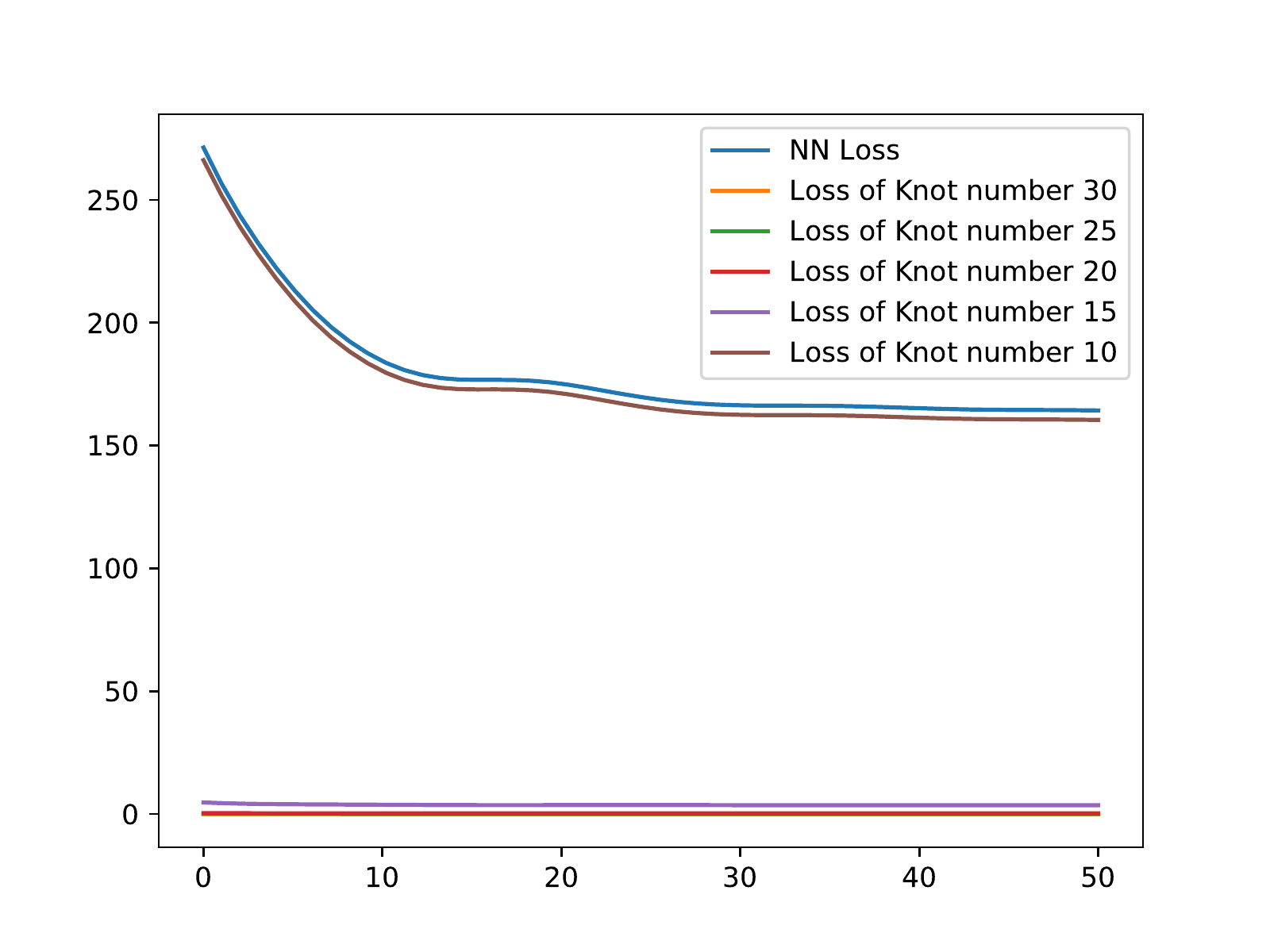}}
\subfigure[Loss curves in the solving phase]{\includegraphics[scale=0.3]{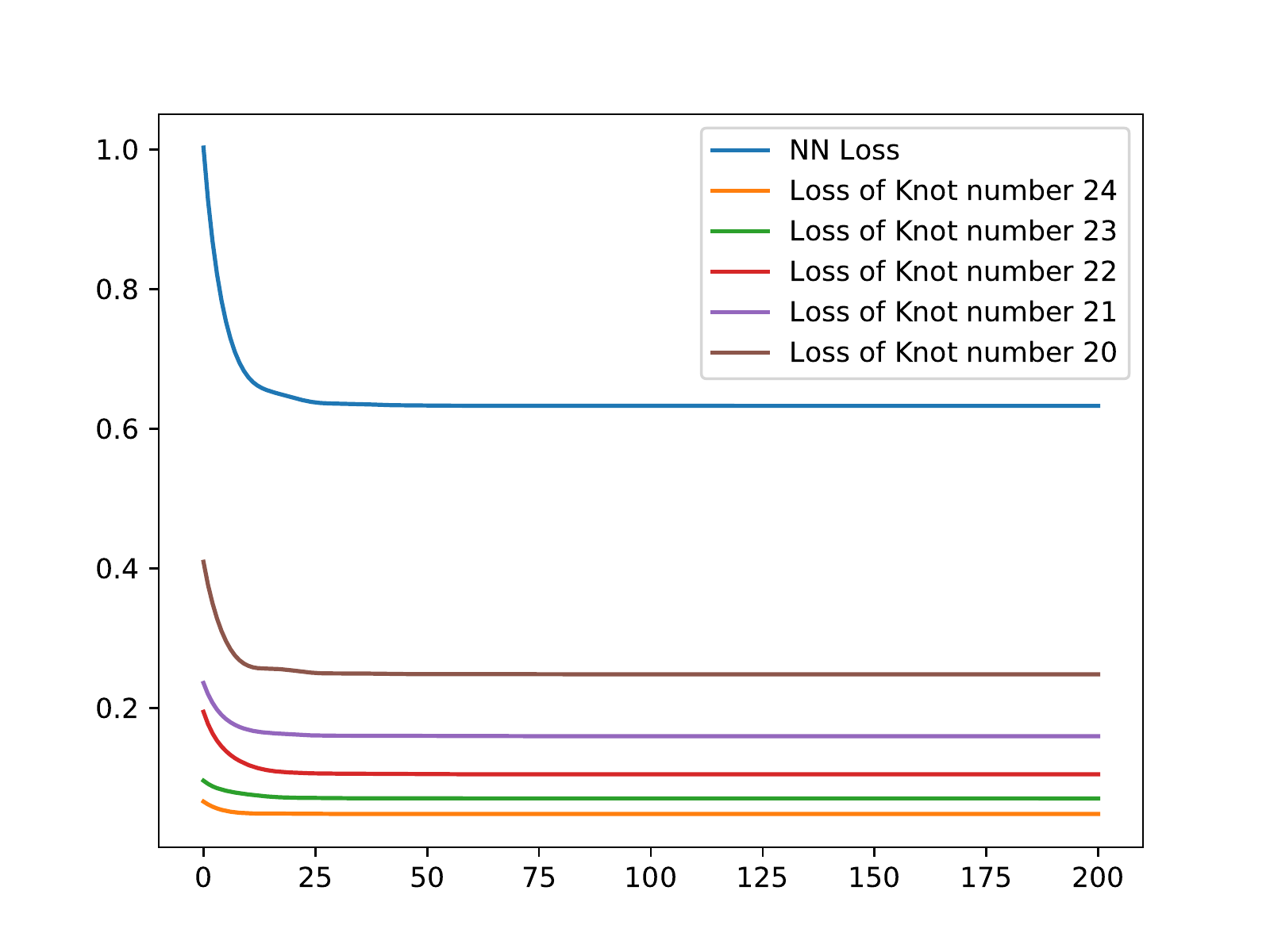}}
\caption{\label{F:demo1} The illustration of how to choose the optimal knot vector based on the loss curves.} 
\end{figure*}

We take Fig.~\ref{F:demo1} as an example to give a further explanation. 
Fig.~\ref{F:demo1}(b) shows the loss functions in the rough train phase for the data given in  Fig.~\ref{F:demo1}(a). We set the range of knot number as $[10,30]$ and $C=5$, that is $dim(U_0) = 30$, $dim(DU_{i+1}) = 30-i\cdot C, i=1,2,3,4$. The loss with $10$ knots is larger than others, this means 10 knots are not sufficient for fitting. In the meanwhile, the loss with $30$ knots is similar to that with $25$ knots, so we choose the knot number range as $[20,25]$ to get a compact lower and upper bound. 
    
Based on the bounds estimated in the rough train, 
$DU_0 = 25$ and $DU_i= 26-i, i=1,2,3,4,5$ are set in the solving phase. There are 5 subnetworks 
($n=5$). Fig.~\ref{F:demo1} (b) shows the loss curves of 5 subnetworks and the whole network. The knot vector $U_4$ which contains 23 interior knots is chosen as the optimal knot vector since the corresponding loss is smaller than that of 20, 21, 22 knots. Here, we do not choose the result of 24 knots since the difference between the two is about the same.

In particular, we emphasize here that the  proposed network in this paper does not need to perform the standard train-predict process.
The network is a solution solver for the problem \eqref{E:topprob} and \eqref{E:inv1}. Once the network is trained by given discrete points $P$, the optimal knot vector is naturally dertermined since the output of the subnetwork is exactly the knot vector.
For different discrete points, the network neededs to be retrained. This differs the standard deep learning methods.

The fitting performance depends on the number of hidden layers in subnetworks. Through lots of numerical experiments, we found that more layers in subnetworks would not significantly improve the performance, but also increase the training time. Hence, we use a network with two hidden layers as a default in this paper.

\begin{remark}
In this paper, we  use cubic B-splines as an illustration to explain the algorithm. The proposed method can be extended to  B-spline curves or surfaces of arbitrary dimension. One only need to change the expression of $A_i$ in Fig.~\ref{F:structure}.
\end{remark}

\subsection{Network architecture integrated parametrization process}
As is well known, the parameterization of discrete points plays a critical role in B-spline curve approximation \cite{ma1995parameterization}. We also adopt a deep neural network to approximate the parametrization based on the same idea of knot  placement problem. 

\begin{figure*}[!htbp]
	\centering
	\subfigure{\includegraphics[scale=0.45]{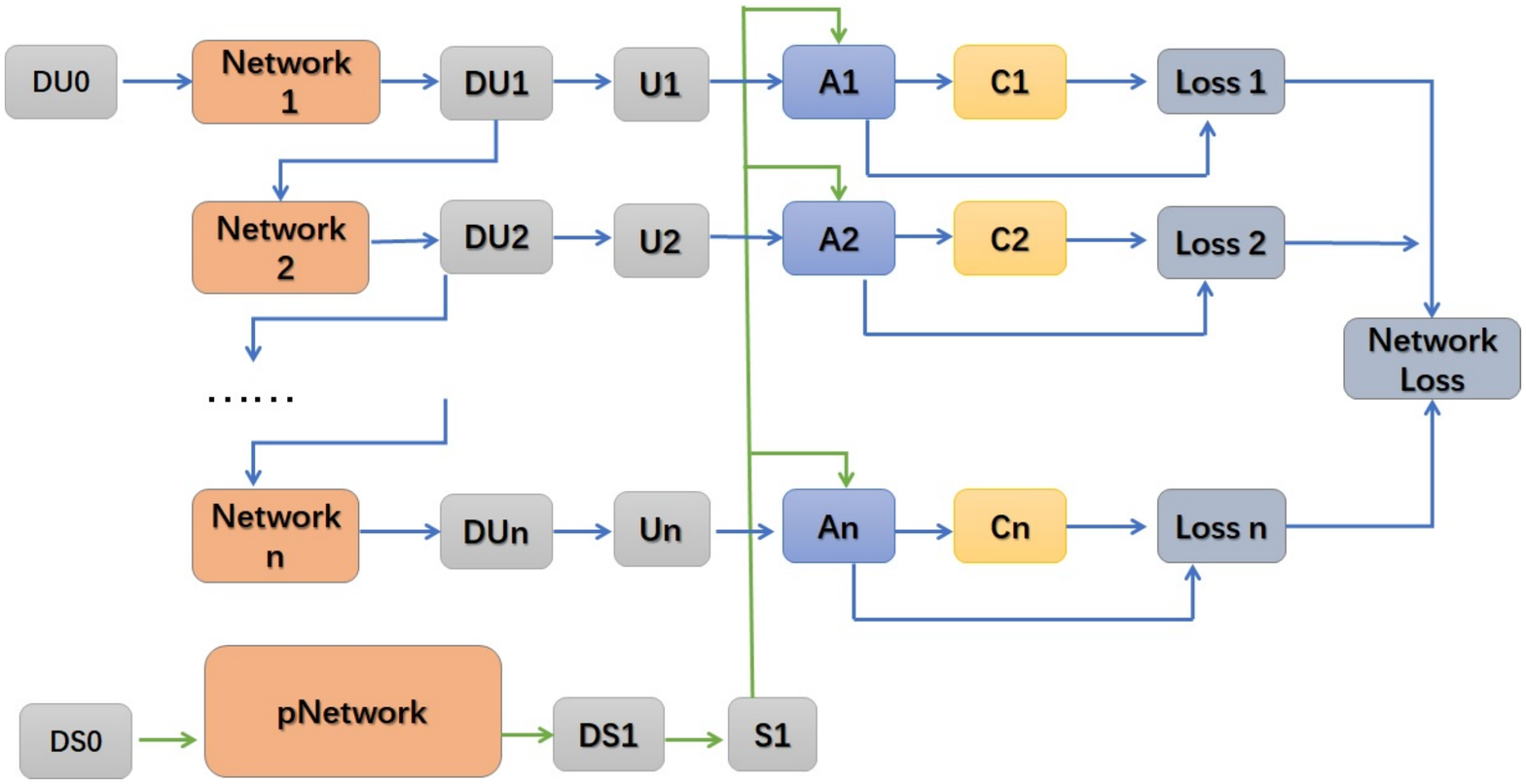}}
	\caption{\label{F:paramstructure} Network architecture integrated parameterization process.}
\end{figure*}

The $pNetwork$ in Fig.~\ref{F:paramstructure} is also a feedforward neural network, which is designed for approximating the optimal parameter values of given data. The input $DS_0$ denotes the parameter difference knot vector, that is $DS_0=\{s_1-s_0, s_2-s_1,\cdots, s_{N}-s_{N-1}\}$.  The output $DS_1$ denotes the difference of the optimal parameter vector $S_1$. The output parameter values $S_1$ is shared by all the subnetwork in computing  the matrix $A_i$. 

\begin{figure*}[!htbp]
	\centering
	\subfigure[Loss curves]{\includegraphics[scale=0.4]{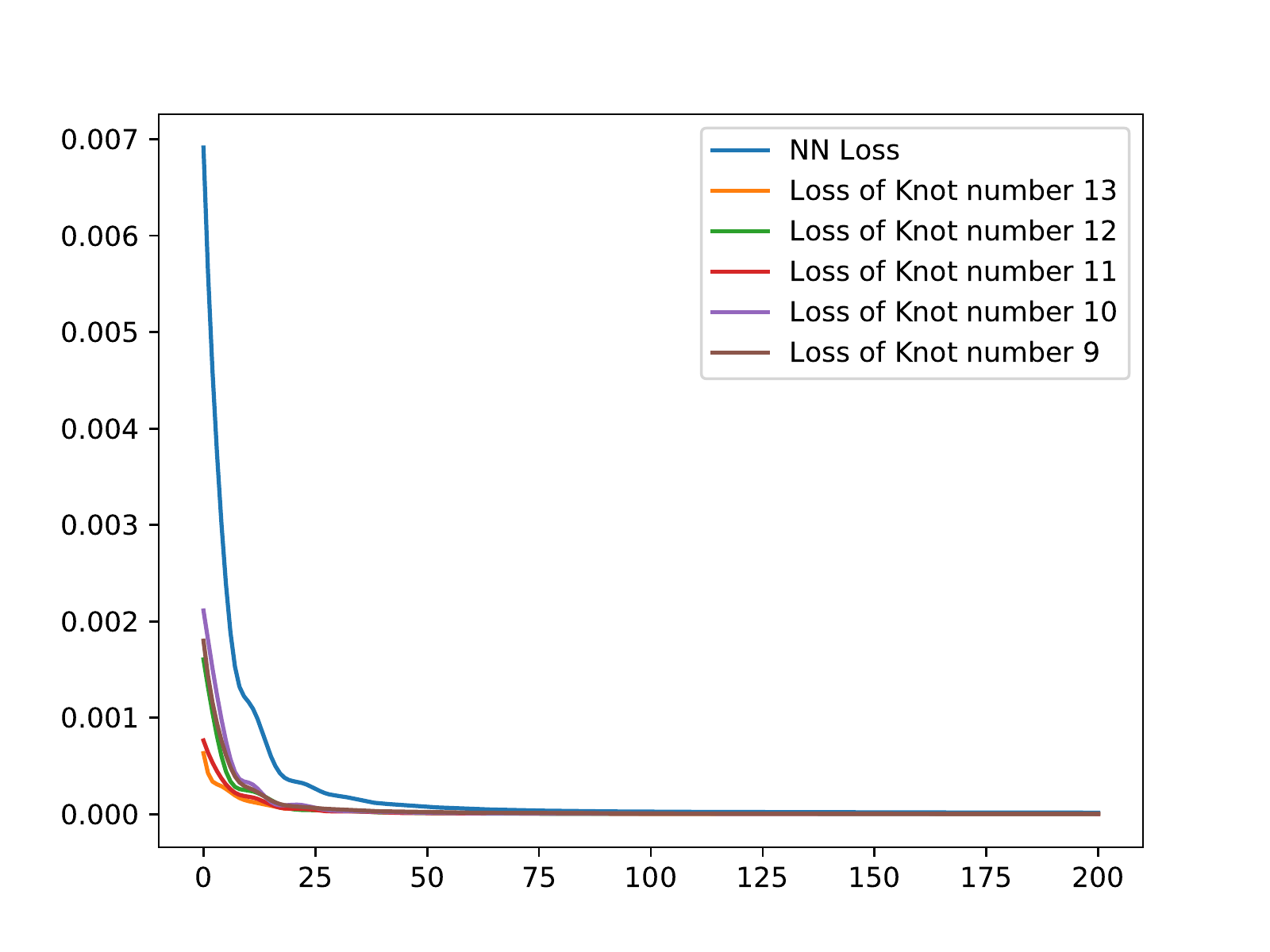}}
	\subfigure[Approximation result with new network]{\includegraphics[scale=0.4]{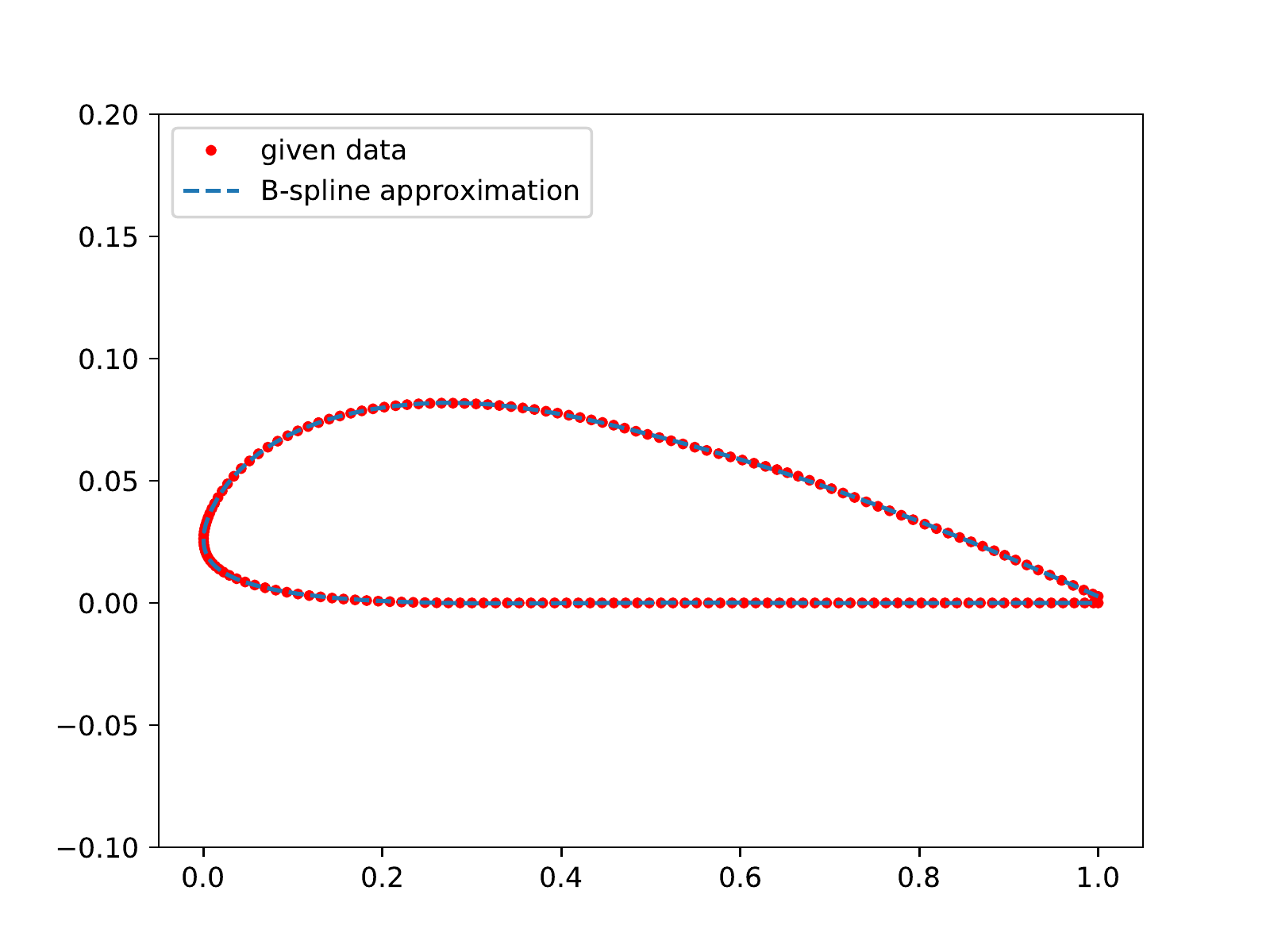}}
	\caption{\label{F:param} B-spline approximation  with the new neural network.}
\end{figure*}

We consider the B-spline curve approximations of the airfoil data considered in Fig.~\ref{F:planar}.
Fig.~\ref{F:param} (a) shows the loss curves obtained by the new network when knot number is chosen as $9-13$.
Fig.~\ref{F:param} (b) shows the B-spline curve approximation with knot number $9$ and fitting error $5.23e-4$.
Compared with the results in Fig.~\ref{F:planar},  the fitting error obtained by the new network improves a little and converges fast, but the time cost (15s) is more because the architecture of the new network is more complicated. In this sense, traditional methods of parameterization should be given priority.

\section{Numerical experiments \label{s.numexp}}
In this section, we present some numerical results to demonstrate the effectiveness of the proposed algorithm. The hausdorff distance is used to measure the fitting performance, which is defined as
$$Dis_H(A, B) = \max\{\sup\limits_{a\in A}\{ Dis(a, B)\},  \sup\limits_{b\in B}\{ Dis(b, A)\}\},$$ 
where $A$ and $B$ are two sets, $Dis(a, B) = \inf \limits_{b\in B} \{ Dis(a,b)\}$, and $Dis(a,b) = \|a-b\|_2$.

The network architecture is simple and it can be trained fast. In our experiments, we use  PyTorch~\cite{paszke2019pytorch,paszke2017automatic} to implement the algorithm. All the  numerical experiments are run on a PC with a 2.3GHz i5-6300HQ cpu and 12GB RAM.

\subsection{Data with noise}
In this section, we consider the B-spline approximations of the function 
$$f(x) = sin(4x-2)+2\exp{^{-30(4x-2)^2}}.$$
We sample 1001 data points uniformly from $f(x)$ and add the noise with uniform distribution in $[-0.05, 0.05]$.
The sparse model introduced in \cite{luo2019knot} and the Lasso model introduced in \cite{yuan2013adaptive} are employed make a comparison.

\begin{figure*}[!htbp]
	\centering
	\subfigure[Loss curves]{\includegraphics[scale=0.4]{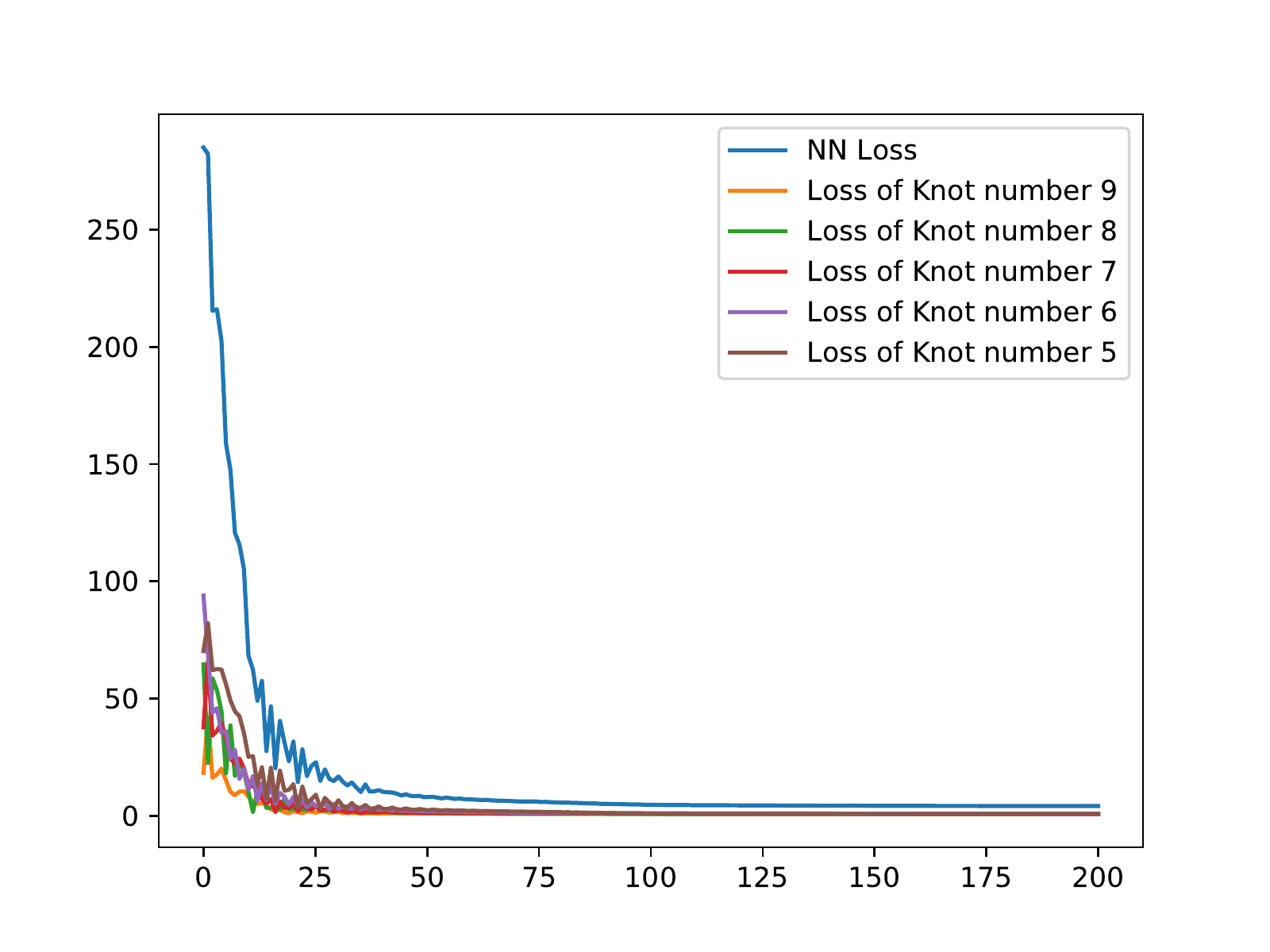}}
	\subfigure[DNN method with $Er=4.9e-2$ and $Kn=5$]{\includegraphics[scale=0.4]{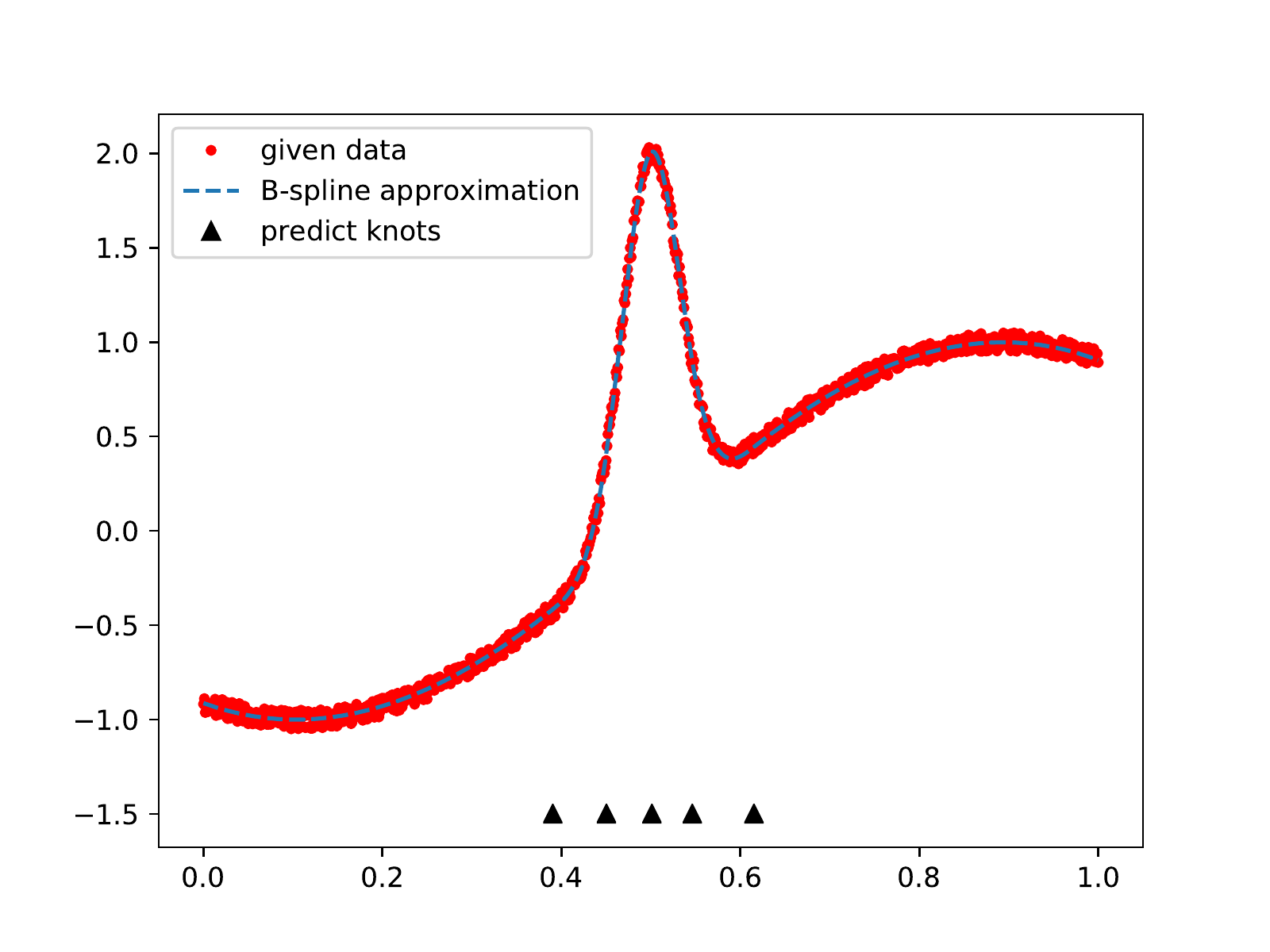}}
	\subfigure[Sparse method with $Er=5.3e-2$ and $Kn=9$]{\includegraphics[scale=0.4]{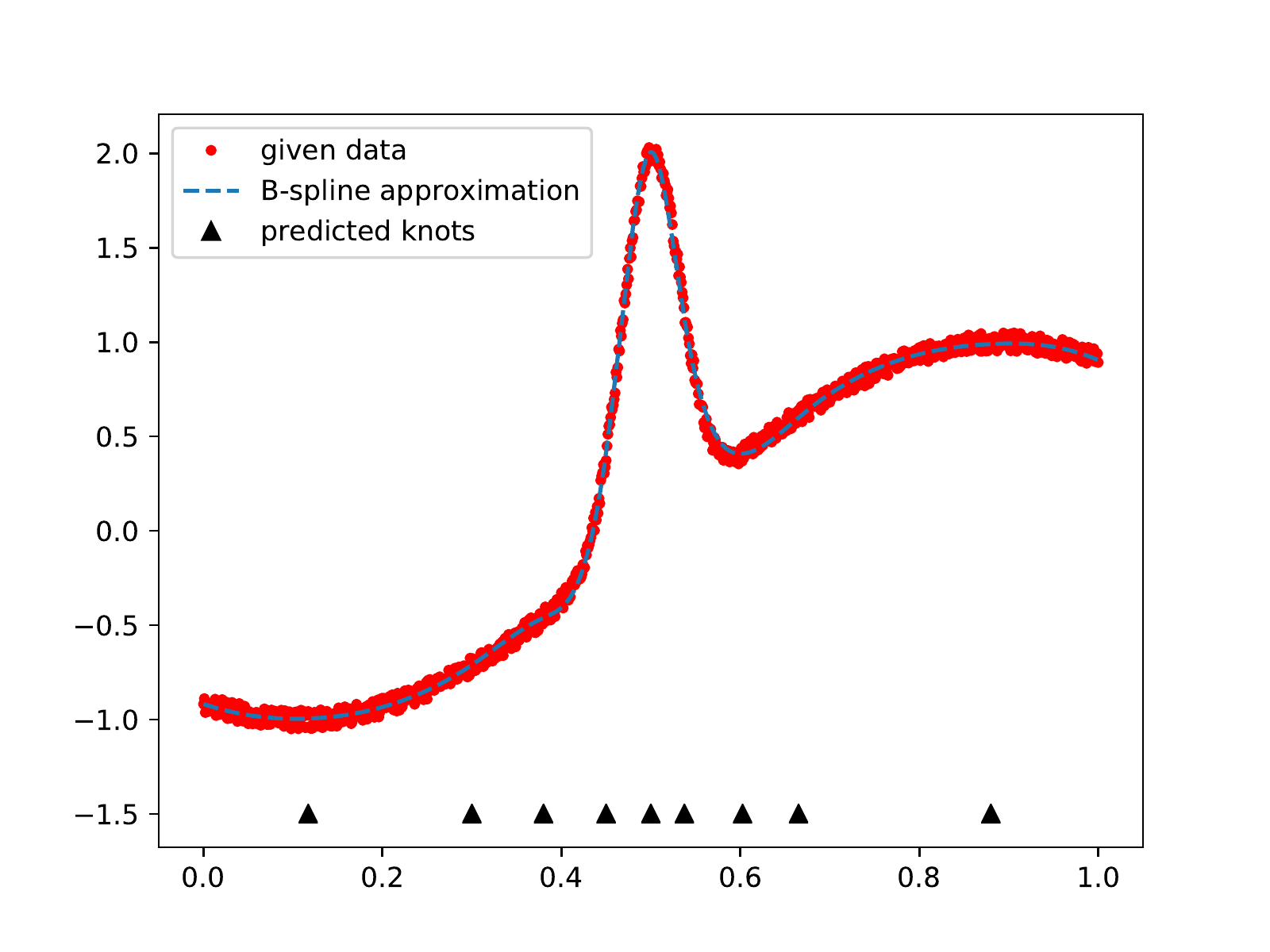}}
	\subfigure[Lasso method with $Er=8.2e-2$ and $Kn=7$]{\includegraphics[scale=0.4]{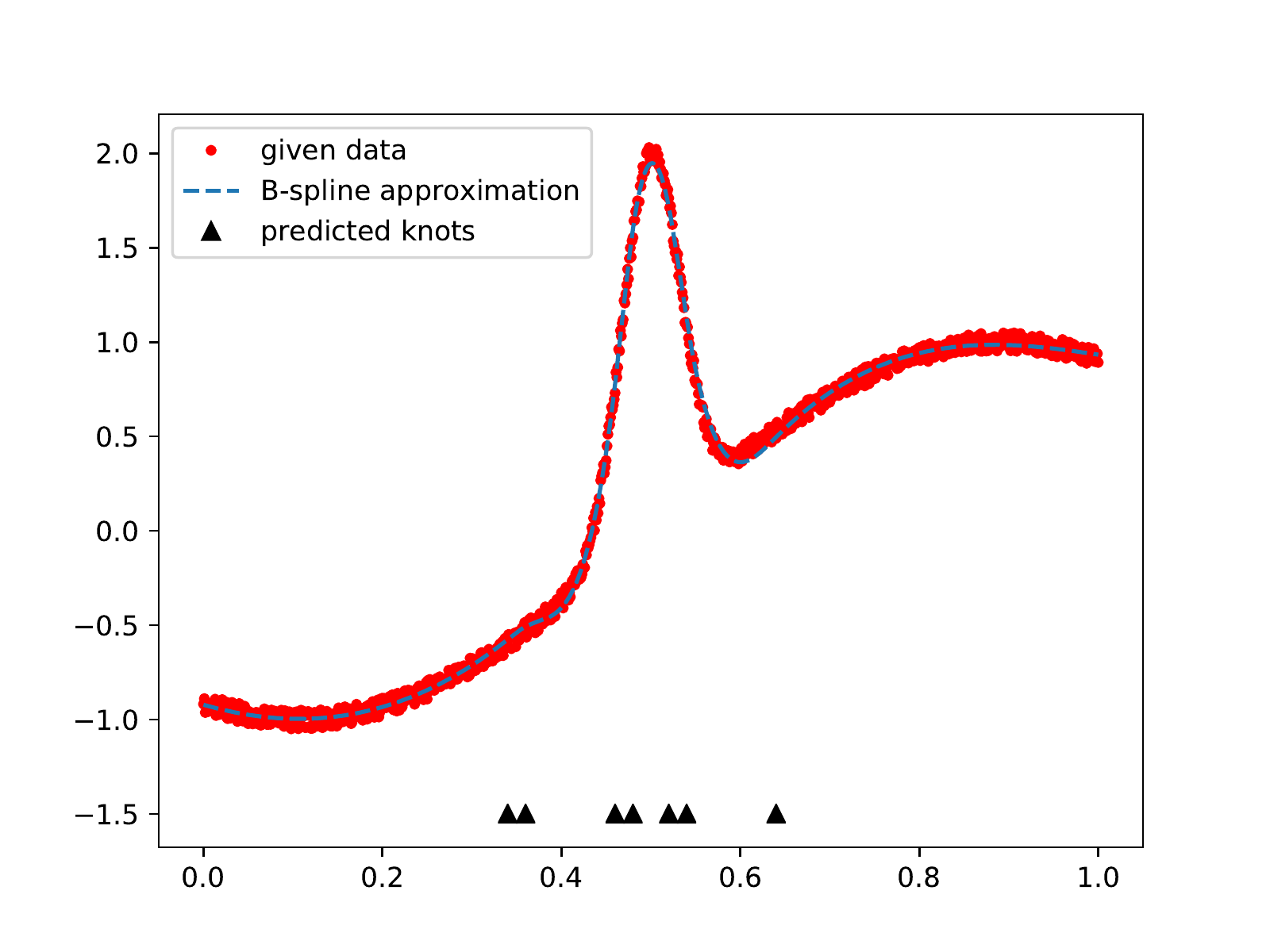}}	
	\caption{\label{F:function} B-spline approximations solved from (b) the proposed method, (c) the sparse model introduced in \cite{luo2019knot}, and (d) the Lasso model \cite{yuan2013adaptive}.}
\end{figure*}

Fig.~\ref{F:function} (a) shows loss curves under different knot number marked by different colors, where the $x-$axis represents the iteration number and $y-$axis represents the loss error of networks. The NN Loss denotes the loss error of the whole neural network in Fig.~\ref{F:function} (a). Since the loss curves of knot number $5,6,7,8,9$ converge to the similar error, thus the minimal knot number $5$ is chosen as the result.
Fig.~\ref{F:function} (b) shows the B-spline approximation with $5$ interior knots computed by the proposed deep neural network method (DNN for short).
Fig.~\ref{F:function} (c) and (d) shows the B-spline approximations solved from the sparse model and Lasso model, respectively.
The underlying knots of B-spline approximations are marked by triangles in Fig.~\ref{F:function} (b-d).

Table.~\ref{T:function} lists the fitting error, time cost and knot number of three methods.
Comparing with the sparse model and Lasso model, our method produces a better approximation with fewer knots and fewer time cost.

\begin{table}[!htbp]
	\label{T:function}
	\caption{Comparisons of approximation performance for a function.}
	\centering
	\begin{tabular}{l|c|c|c}	
		\hline
		methods  & error (Er)  & time & knot number (Kn)   \\
		\hline
		Sparse model &5.3e-2  & 22s & 9 \\
		\hline
		Lasso model & 8.2e-2  & 180s & 7 \\
		\hline
		Our method & 4.9e-2& 8s & 5\\
		\hline
	\end{tabular}
\end{table}

\subsection{Planar curve}
\begin{figure*}[!htbp]
	\centering
	\subfigure[Loss curves]{\includegraphics[scale=0.4]{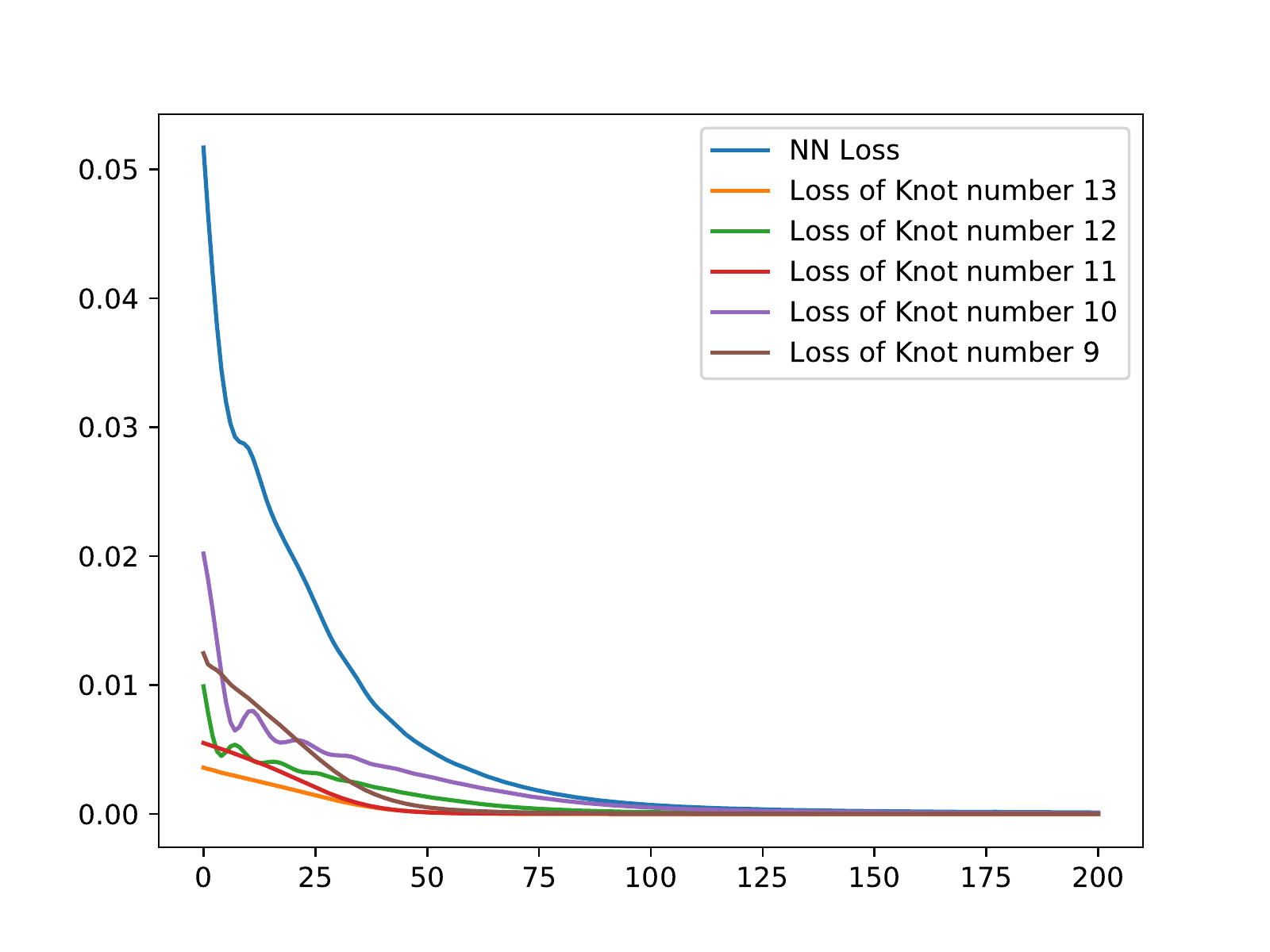}}
	\subfigure[DNN method]{\includegraphics[scale=0.4]{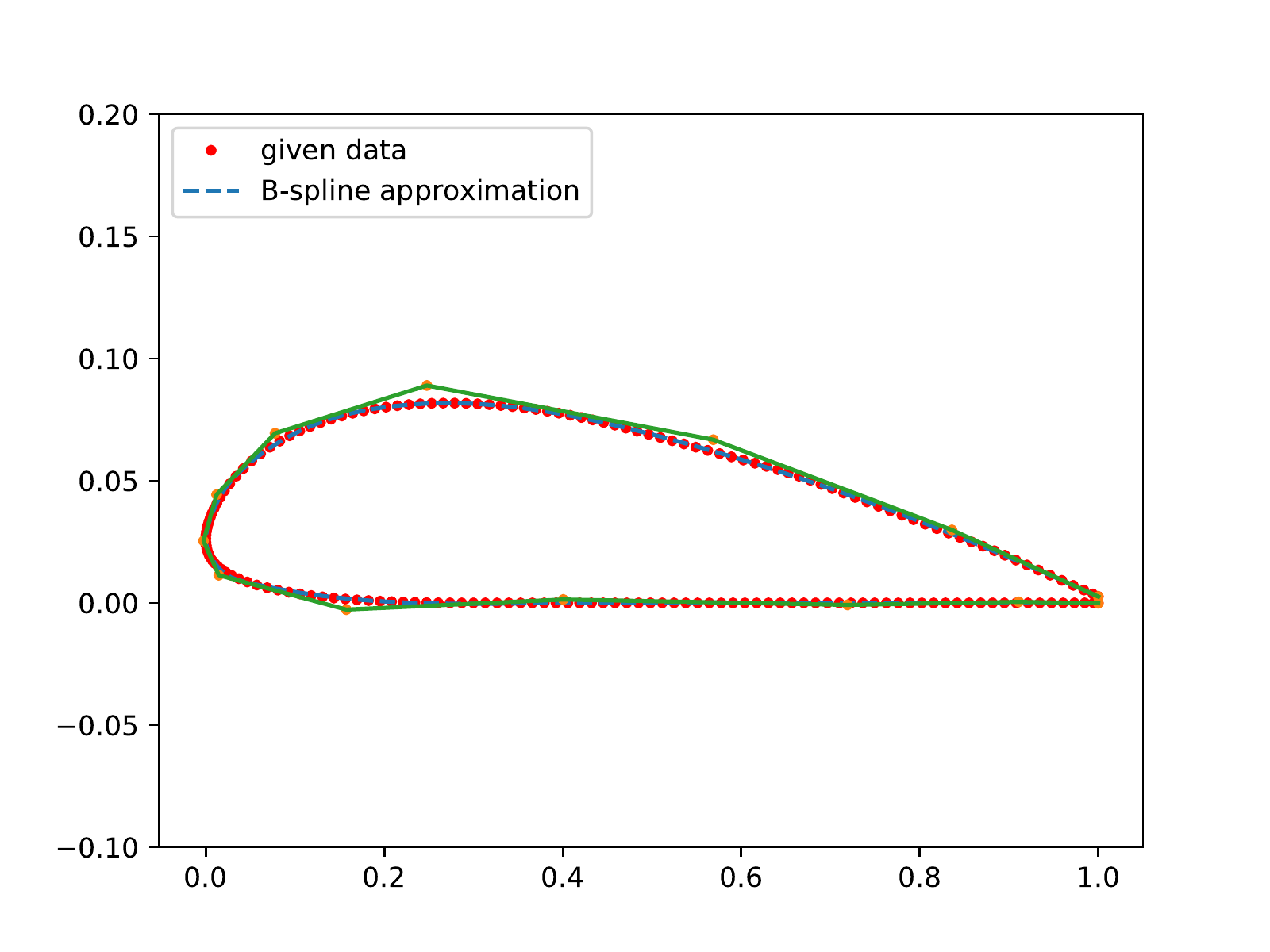}}
	\subfigure[DOM method]{\includegraphics[scale=0.4]{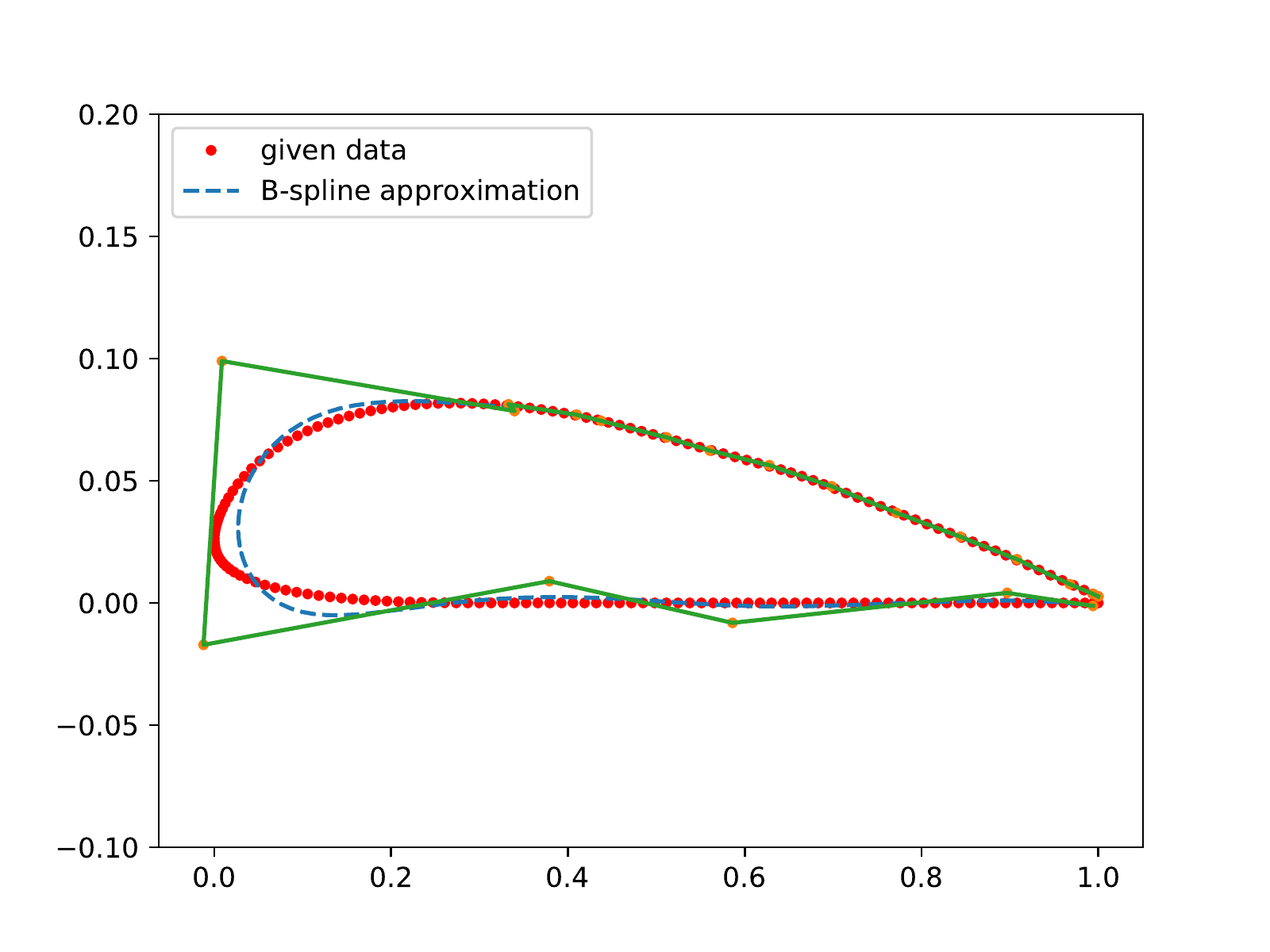}}
	\subfigure[NKTP method]{\includegraphics[scale=0.4]{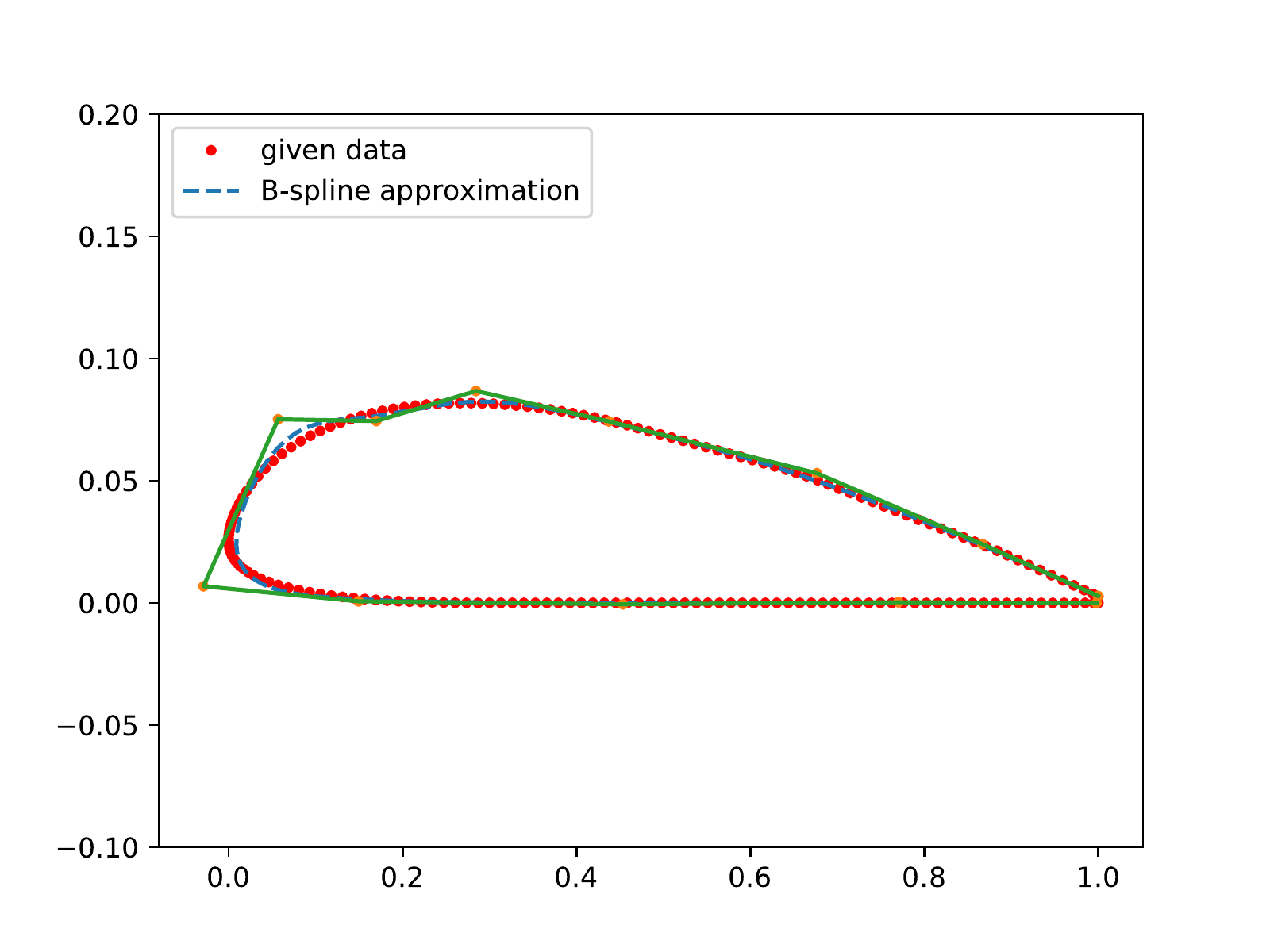}}
	\caption{\label{F:planar} B-spline curve approximation of the airfoil data constructed from different methods.}
\end{figure*}

In this subsection, an airfoil data is considered to verify the practical ability of the deep neural network based method. 
We firstly use accumulated chord length parameterization to parameterize the data. 
Fig.~\ref{F:planar}(a) shows the loss curves under different knot number. Since the loss curves of knot number $9-13$ converge to the similar error, thus the minimal knot number $9$ is chosen as the result. Fig.~\ref{F:planar}(b) shows the B-spline approximation with $9$ interior knots and fitting error $6.7e-4$ computed by the deep learning based method.

The DOM method \cite{park2007b-spline} and NKTP method \cite{piegl2000surface} are compared with our method. The DOM method needs a pre-approximation in order to calculate curvature thus the fitting performance depends on the pre-approximation. The NKTP method calculates knots based on parameter values of given data and needs a pre-specified knot number. The statistics of DOM method, NKTP method and our method are summarised in Table ~\ref{T:planar}. It can be seen that although the other two methods are very fast, the posterior fitting error is not considered in the entire calculation process, so the fitting performance is not as good as ours.

\begin{table}[!htbp]
	\label{T:planar}
	\caption{Comparisons of approximation performance of the airfoil data.}
	\centering
	\begin{tabular}{l|c|c|c}	
		\hline
		methods  & error (Er)  & time & knot number  (Kn) \\
		\hline
		DOM & 2.7e-2  &0.5s & 15\\
		\hline
		NKTP&  6.4e-3  &5e-4s &	10 \\
		\hline
		Our method & 6.7e-4 &5s &9\\
		\hline
	\end{tabular}
\end{table}

\subsection{Space curve}
We consider the B-spline curve approximation of the $C^0$ continuous curve in $\mathbb{R}^3$ defined by
\begin{eqnarray*}
	x(t) &=& \frac{100}{\exp{^{|10t-5|}}}+\frac{25(2t-1)^5}{4},\\
	y(t) &=& 30t-5,\\
	z(t) &=& 100t^2, ~~~ t\in[0, 1].
\end{eqnarray*}

\begin{figure*}[!htbp]
	\centering
	\subfigure[Loss curves]{\includegraphics[scale=0.4]{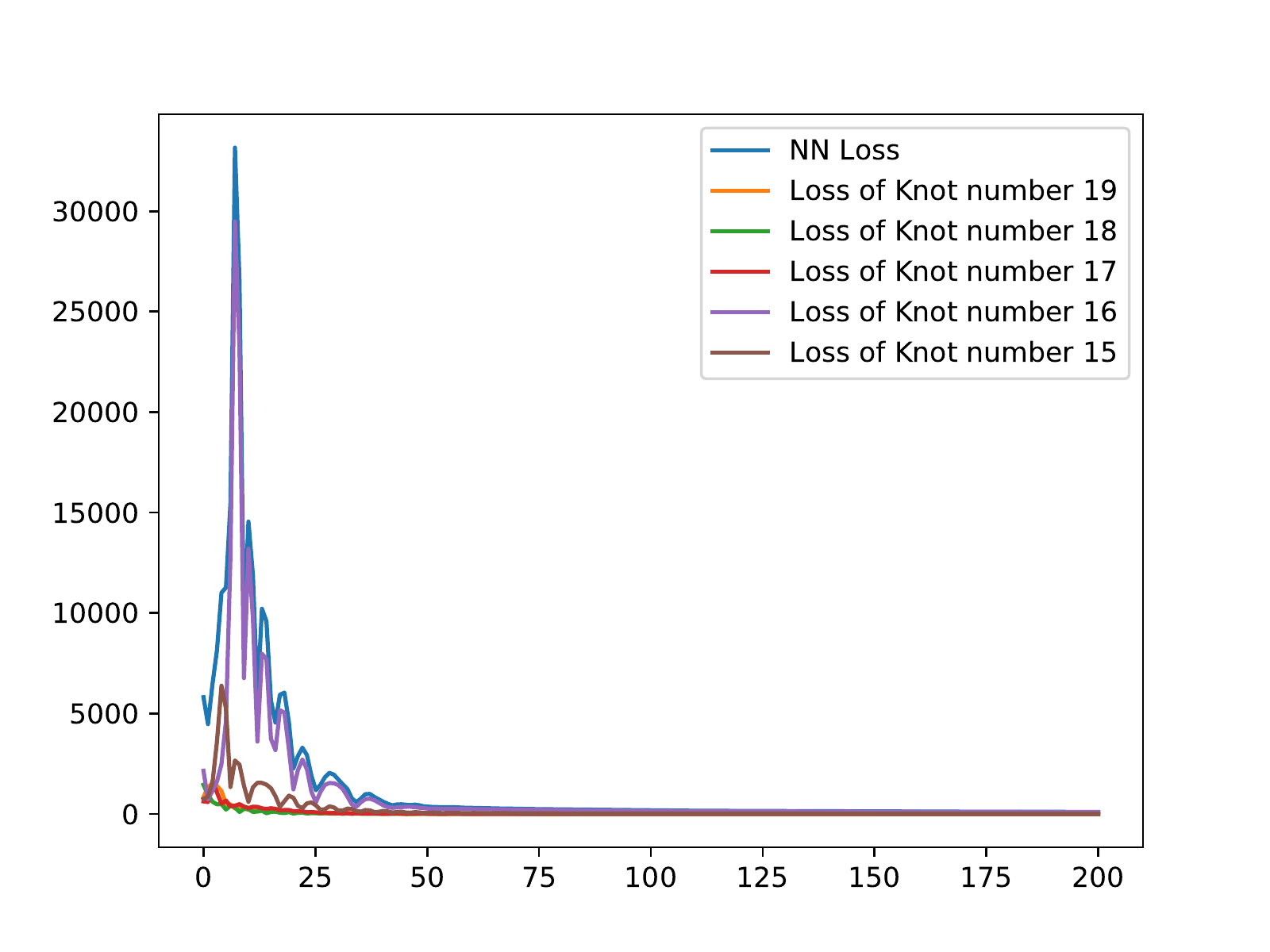}}
	\subfigure[B-spline curve approximation ]{\includegraphics[scale=0.4]{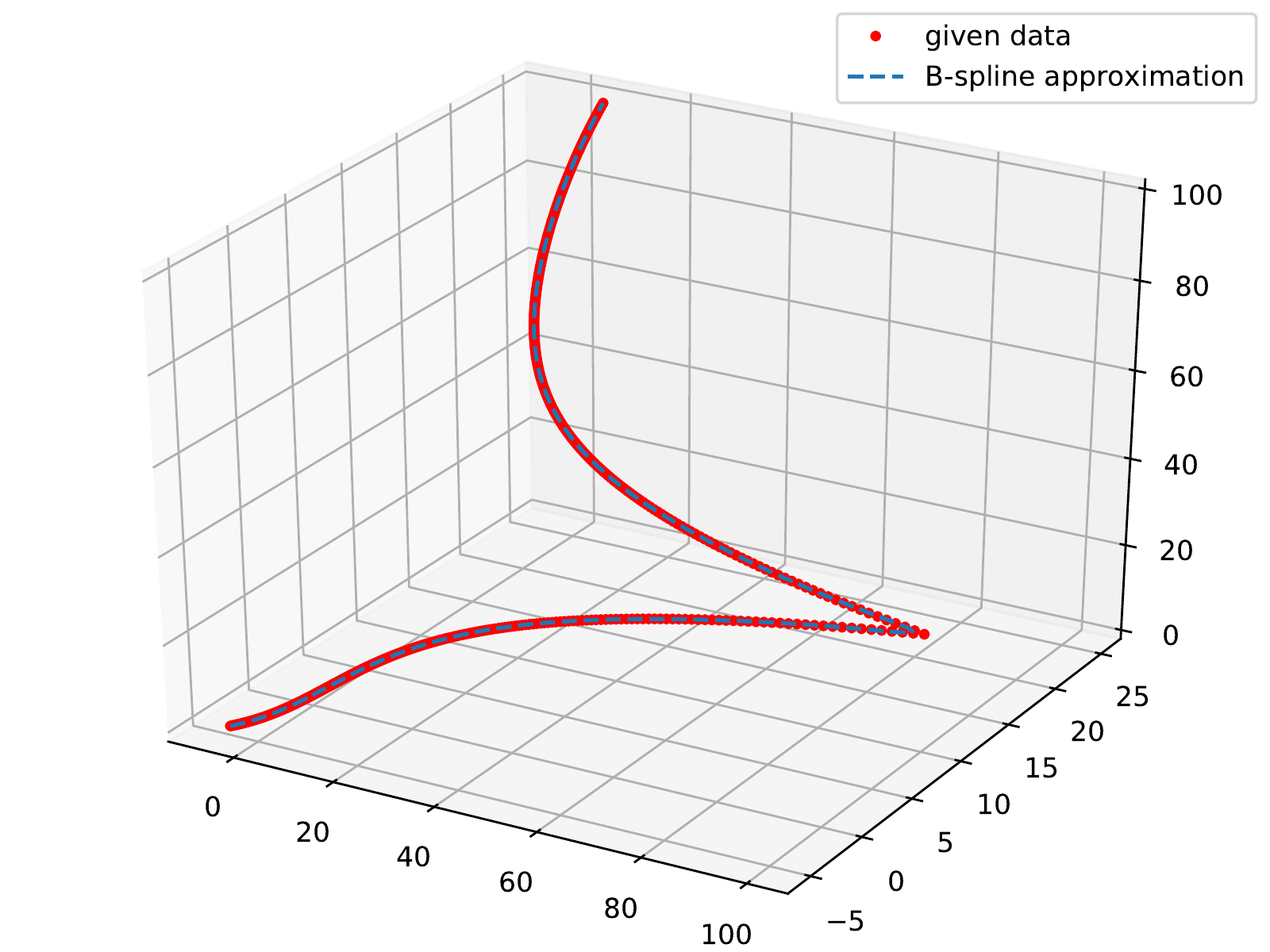}}
	
	\caption{\label{F:space} The B-spline curve approximation of a space curve and the correpsonding loss curves.}
\end{figure*}

We make a comparison with the knot removal method \cite{lyche1987knot} and the genetic algorithm \cite{yoshimoto2003data}. There are $501$ data points which are sampled uniformly from the above curve. According to the loss curve shown in Fig.~\ref{F:space} (a), the minimal knot number $15$ and the corresponding knots are chosen as the result. The resulting B-spline curve approximation with fitting error $Er=0.17$ is shown in Fig.~\ref{F:space} (b). Table~\ref{T:space} summarizes the statistics of different methods. Compared with knot removal method and genetic algorithm, our method has an overwhelming advantage in knot number and fitting performance. Besides, it is much faster than other two methods, this is because the compared methods need to compute the least squares fitting error in each iteration, which is very time-consuming.

Knot removal has the advantages of finding the accurate multiple knots, but the process is artificially tried, while other two methods are automatic. The limitation of our DNN method is that it can not find the accurate multiple knots since the Softmax function never equal to 0, so if we want to get multiple knots, we need to post-process the results. 

\begin{table}[!htbp]
	\label{T:space}
	\caption{Comparisons of approximation performance of a space curve.}
	\centering
	\begin{tabular}{l|c|c|c|c}	
		\hline
		methods  & error (Er) & time & knot number (Kn)& multiple knots \\
		\hline
		Knot removal & 0.15  & 624s & 24 & 0.5,0.5,0.5\\
		\hline
		Genetic algorithm  & 0.42  & 182s &	 41 & 0.499,0.505,0.507\\
		\hline
		Our method & 0.17 & 8s & 15 & 0.489, 0.501, 0.509\\
		\hline
	\end{tabular}
\end{table}

\section{Conclusion \label{s.conclusion}}
In this paper, we propose a deep neural network framework for B-spline curve approximation. Different from traditional methods based on optimization or local features, the knot placement here is understood as an approximation process of initial knots. It is a bit tricky to handle the approximation by traditional techniques due to the nonlinearity and non-convexity between knots and B-spline curves. Deep neural network is proven to has powerful approximation capabilities. Thus, we design a neural network composed of several subnetworks to approximate input initial knots. Each subnetwork outputs the optimal knot positions with fixed knot number. Then, all the subnetworks stacked together approximate knots of different dimensions. The number of subnetworks can be specified by users to satisfy fitting error. In addition, the parameterization of given data can be integrated into the neural network easily.

For the future work, the similar idea can be extended in the parameterization problem in isogemetric analysis field \cite{hughes2005isogeometric}.
Computing a valid parameterization of the geometric model represented by the boundary is a fundamental and crucial problem in isogemetric analysis. A deep neural network designed suitably may be helpful for addressing this tough problem.


\section*{References}
\bibliographystyle{unsrt}
\bibliography{DNN_knot}

\end{document}